\documentclass{ifacconf}
\usepackage{graphicx}
	\graphicspath{{./}{figure/}}
\usepackage{natbib}
\usepackage{color}
\usepackage{amsmath,amsfonts,amssymb,paralist}

\DeclareMathOperator*{\argmin}{arg\,min}
\newcommand{\abs}[1]{\left\vert#1\right\vert}
\newcommand{\cE}{\mathcal{E}}
\newcommand{\cU}{\mathcal{U}}
\newcommand{\cV}{\mathcal{V}}
\newcommand{\R}{\mathbb{R}}

\begin{document}
\begin{frontmatter}

\title{Control strategies for road risk mitigation in kinetic traffic modelling}

\author[First]{A. Tosin}
\author[Second]{M. Zanella} 

\address[First]{Department of Mathematical Sciences ``G. L. Lagrange'' \\
	Politecnico di Torino, Torino, Italy \\
	(e-mail: andrea.tosin@polito.it)}
\address[Second]{Department of Mathematical Sciences ``G. L. Lagrange'' \\
	Politecnico di Torino, Torino, Italy \\
	(e-mail: mattia.zanella@polito.it)}

\begin{abstract}
In this paper we present a Boltzmann-type kinetic approach to the modelling of road traffic, which includes control strategies at the level of microscopic binary interactions aimed at the mitigation of speed-dependent road risk factors. Such a description is meant to mimic a system of driver-assist vehicles, which by responding locally to the actions of their drivers can impact on the large-scale traffic dynamics, including those related to the collective road risk and safety.
\end{abstract}

\begin{keyword}
Road traffic, traffic control, binary interactions, Boltzmann-type equations
\end{keyword}

\end{frontmatter}

\section{Introduction}
Our inner city mobility is rapidly changing due to the automation of driving and the sharing of information and communication technology. This process is leading to the creation of new paradigms in terms of efficient infrastructure and traffic management solutions. Among others, we mention in this direction the broad developments in the technology for driver-assist cars, self-driving cars and intelligent intersections, see~\cite{santi2014PNAS,tachet2016PLOS1}. As an effect of the fast rising of urban population, such an automation process also shed light on safety issues in road traffic management. According to recent reports on traffic safety in the world, see e.g.~\cite{WHO2004report,WHO2015report}, road risk arises as a result of several factors largely linked to the subjectivity of the driving behaviour of the individuals. Among others, here we recall in particular those related to the variability of the speed in the traffic flow: large differences in the speeds of the vehicles within the traffic stream are reported to be responsible for an increase in the crash risk.

So far, road risk and safety have been mainly investigated by means of empirical approaches. These include, for instance, the analysis of the distribution of the fatality rates over time or the study of the accident time series and of safety indicators, see e.g.~\cite{oppe1989AAP,hermans2008AAP,hermans2009AAP}. Nevertheless, recently theoretical efforts have been devoted to the comprehension of the links between traffic dynamics and safety issues by means of mathematical models, see e.g.~\cite{herty2011ZAMM,moutari2013IMAJAM,moutari2014CMS,freguglia2017CMS}.

In this paper we continue along the latter research line by combining a Boltzmann-type kinetic description of the road traffic, cf.~\cite{klar1997JSP,herty2010KRM,puppo2016CMS} for related approaches, with a preliminary study of control strategies in the frame of the driver-assist car technology for the mitigation of the risk caused by the speed variance of the vehicles. The kinetic approach is particularly appropriate to our goal thanks to its fundamental link with the particle representation of the driver-vehicle system, which is precisely the level at which driver-assist control strategies can act. At the same time, it allows one to upscale rigorously such small-scale dynamics to an aggregate level, which is more suited to engineering needs.

The control approach adopted here has roots in the Model Predictive Control (MPC), which has been used in the engineering community since over fifty years, see  e.g.~\cite{camacho2007BOOK,michalska1993IEEE,sontag1998BOOK} for an overview and further references. MPC methods have been traditionally employed in the frame of ODEs, whereas for kinetic and fluid dynamic equations few results are available in the literature, cf.~\cite{albi2014PTRSA,albi2015CMS}. The hallmark of the kinetic formulation of the control problem is the derivation of an explicit feedback control for binary, i.e. one-to-one, vehicle dynamics, which is then straightforwardly embedded into a Boltzmann-type kinetic equation for a large number of vehicles. As it is well known, MPC leads typically to a control which is suboptimal with respect to the theoretical optimal one. Nevertheless, performance bounds can be established which guarantee the consistency of such an approximation in the kinetic framework, see~\cite{grune2009SICON,herty2017DCDS}. In addition to that, the proposed Boltzmann formulation of the MPC has an overall computational cost which scales linearly with the total number of vehicles of the system. This makes it competitive compared to other techniques for computing the optimal control.

In more detail, the paper is organised as follows. In Section~\ref{sect:binary} we present the unconstrained microscopic traffic dynamics via the concept of binary interactions. In Section~\ref{sect:control} we introduce the binary control and discuss possible strategies for speed-dependent road risk mitigation. In Section~\ref{sect:boltzmann} we embed the constrained microscopic dynamics into a kinetic Boltzmann-type equation, which we then use to investigate analytically the large-scale impact of the envisaged risk mitigation strategies. In Section~\ref{sect:numerics} we provide numerical evidences of the risk mitigation effect by simulating the fundamental diagrams of traffic with special focus on the evolution of the speed variance. Finally, in Section~\ref{sect:conclusion} we summarise the main aspects of the proposed approach and we briefly sketch research perspectives.

\section{Microscopic binary interactions}
\label{sect:binary}
The kinetic modelling approach relies on the concept of \emph{binary interactions} at the particle level, which fits naturally the \emph{follow-the-leader} principle that most microscopic models of vehicular traffic are based on, cf.~\cite{gazis1961OR}.

We describe the microscopic state of a vehicle by a scalar variable $v\in [0,\,1]$ representing the (dimensionless) speed. If $w\in [0,\,1]$ is the speed of the leading vehicle, we assume that in a short time interval $\Delta{t}>0$ an interaction between the two vehicles produces a change of speed of the former described by the rule
\begin{equation}
	v'=v+\Delta{t}I(v,\,w;\,\rho),
	\label{eq:binary}
\end{equation}
where $v'$ is the post-interaction speed and
\begin{equation}
	I(v,\,w;\,\rho):=
	\left\{
	\begin{array}{ll}
		P(\rho)\left(\min\{v+\Delta{v},\,1\}-v\right) & \textnormal{if } v<w \\[2mm]
		(1-P(\rho))\left(P(\rho)w-v\right) & \textnormal{if } v>w
	\end{array}
	\right.
	\label{eq:I}
\end{equation}
is the interaction function. In particular, $\Delta{v}>0$ is the increase in speed when the vehicle accelerates, $\rho\in [0,\,1]$ is the (dimensionless) macroscopic density of the vehicles and
\begin{equation}
	P(\rho):=1-\rho^\gamma, \quad \gamma>0,
	\label{eq:P}
\end{equation}
is the probability of accelerating. The function~\eqref{eq:I} expresses the fact that a vehicle accelerates if it is slower than the leading vehicle ($v<w$) and brakes if it is faster ($v>w$). In the former case it increases its speed by a quantity which is at most $\Delta{v}$ (the ``$\min$'' guarantees that the bound $v'\leq 1$ is preserved) while in the latter case it decreases its speed to a fraction $P(\rho)w$ of the speed of the leading vehicle. Owing to~\eqref{eq:P}, the lighter the traffic (i.e. low $\rho$) the closer to $w$ the speed targeted when breaking. Finally, acceleration and breaking are more or less probable depending on the congestion of the traffic, which is expressed by the coefficients $P(\rho)$ and $1-P(\rho)$ in~\eqref{eq:I}.

The leading vehicle is instead assumed not to change speed in consequence of the interaction just described, because binary interactions in vehicular traffic are mainly anisotropic. Therefore we set
$$ w'=w. $$

Notice that the binary interaction rules for $v,\,w$ can be seen as a time discretisation of the following equations:
$$ \frac{dv}{dt}=I(v,\,w;\,\rho), \qquad \frac{dw}{dt}=0 $$
relating the acceleration of a car to the interaction with its leading vehicle in a time interval $(t,\,t+\Delta{t}]$.

\section{Binary control strategies for risk mitigation}
\label{sect:control}
Having in mind driver-assist vehicles, we now include in the previous setting a reaction ability of the cars to the actions of the drivers aimed at enhancing the driving safety. Thus we modify the binary interaction rules set forth in Section~\ref{sect:binary} by adding a control term $u$ such that
\begin{equation}
	\frac{dv}{dt}=I(v,\,w;\,\rho)+u, \qquad \frac{dw}{dt}=0.
	\label{eq:binary_continuous.u}
\end{equation}
The control is supposed to be applied by the car in response to the changes of speed imposed by the driver so as to minimise a certain cost functional $J=J(v,\,w,\,u)$ linked to a measure of the driving risk. Hence the optimal control $u^\ast$ is defined by
\begin{equation}
	u^\ast:=\argmin_{u\in\cU}J(v,\,w,\,u)
	\label{eq:ustar}
\end{equation}
subject to~\eqref{eq:binary_continuous.u}, $\cU$ being a set of admissible controls to be suitably specified.

Since the differences in the speed of the vehicles along the road have been recognised as a non-negligible factor of driving risk, cf.~\cite{WHO2004report}, a conceivable form of the cost functional $J$ to be minimised may be one which involves the \emph{binary variance} of the speeds of the interacting vehicles. This leads us to consider:
\begin{equation}
	J(v,\,w,\,u)=\frac{1}{2}\int_t^{t+\Delta{t}}\left[(w-v)^2+\nu u^2\right]\,ds,
	\label{eq:J.variance}
\end{equation}
where the term $\frac{1}{2}(w-v)^2$ is the aforesaid binary variance while $\frac{\nu}{2}u^2$, $\nu>0$, is a penalisation on large controls.

Another option is to minimise the gap between the current speed of the car and a certain desired (or imposed) speed $v_d\in [0,\,1]$, which may be understood for instance as a speed limit or as a recommended speed fostering the occurrence of green waves. In this case we may consider the cost functional
\begin{equation}
	J(v,\,w,\,u)=\frac{1}{2}\int_t^{t+\Delta{t}}\left[\left(v_d-v\right)^2+\nu u^2\right]\,ds,
	\label{eq:J.vd}
\end{equation}
cf.~\cite{albi2015CMS} in a different context.

\subsection{Feedback control}
In order to tackle the control problem~\eqref{eq:binary_continuous.u}-\eqref{eq:ustar} we should consider a bounded control $-\infty<a\leq u\leq b<+\infty$. The values $a,\,b$ should guarantee that the bounds $0\leq v\leq 1$ on the post-interaction speed resulting from~\eqref{eq:binary_continuous.u} are not violated in the whole time interval $(t,\,t+\Delta{t}]$. However, instead of considering the constrained minimisation problem~\eqref{eq:ustar} we will admit that $u\in\R$ and we will show that it is possible to preserve the aforesaid bounds by carefully selecting $\Delta{t}$ and $\nu$ \emph{a posteriori}.

We consider at first the cost functional~\eqref{eq:J.variance}. The Hamiltonian of the control problem~\eqref{eq:binary_continuous.u}-\eqref{eq:ustar} is in this case:
$$ H(v,\,w,\,u,\,\lambda):=\frac{1}{2}(w-v)^2+\frac{\nu}{2}u^2+\lambda\left(I(v,\,w;\,\rho)+u\right), $$
$\lambda=\lambda(t)$ being the Lagrange multiplier. From Pontryagin's principle, the optimality conditions turn out to be:
\begin{equation*}
	\begin{cases}
		\nu u+\lambda=0 \\[1mm]
		\frac{d\lambda}{dt}=w-v-\lambda\partial_vI(v,\,w;\,\rho) \\[2mm]
		\lambda(t+\Delta{t})=0,
	\end{cases}
\end{equation*}
which we discretise in $(t,\,t+\Delta{t})$ as
\begin{equation*}
	\begin{cases}
		\nu u+\lambda=0 \\
		\lambda'=\lambda+\Delta{t}\left(w'-v'-\lambda'\partial_vI(v',\,w';\,\rho)\right) \\
		\lambda'=0.
	\end{cases}
\end{equation*}
We have denoted by $'$ the variables computed at $t+\Delta{t}$ and, in particular, we have used the implicit Euler scheme for the equation of the multiplier. As a result we get
$$ u=\frac{\Delta{t}}{\nu}(w'-v') $$
where $v',\,w'$ have to be understood as the post-interaction speeds produced by the constrained binary interaction rules resulting from the time discretisation of~\eqref{eq:binary_continuous.u}, i.e.
\begin{align}
	\begin{aligned}[c]
		v' &= v+\Delta{t}I(v,\,w;\,\rho)+\Delta{t}u \\
		w' &= w.
	\end{aligned}
	\label{eq:binary.u}
\end{align}
Using these expressions we deduce
\begin{equation}
	u=\frac{\Delta{t}}{\nu+\Delta{t}^2}(w-v)-\frac{\Delta{t}^2}{\nu+\Delta{t}^2}I(v,\,w;\,\rho),
	\label{eq:u.variance}
\end{equation}
namely we get $u$ in \emph{feedback form} as a function of the pre-interaction speeds. We notice that, consistently with the MPC approach together with a receding horizon strategy, the control $u$ in~\eqref{eq:binary.u} is assumed to be constant in the time horizon $\Delta{t}$ coinciding with the characteristic time of a binary interaction.

By plugging~\eqref{eq:u.variance} into~\eqref{eq:binary.u} we finally deduce the following feedback-constrained binary interaction scheme:
\begin{align}
	\begin{aligned}[c]
		v' &= v+\dfrac{\nu\Delta{t}}{\nu+\Delta{t}^2}I(v,\,w;\,\rho)+\frac{\Delta{t}^2}{\nu+\Delta{t}^2}(w-v) \\
		w' &= w
	\end{aligned}
	\label{eq:binary.u.variance}
\end{align}
corresponding to the instantaneous strategy of reducing the speed variance of the interacting vehicles. Using the expression~\eqref{eq:I} of the interaction function $I$ it is possible to check that if $0<\Delta{t}\leq 1$ then $v'\in [0,\,1]$ for any given $v,\,w\in [0,\,1]$ and any $\nu>0$. In particular, if $\nu\to+\infty$ then~\eqref{eq:binary.u.variance} reduces to the unconstrained binary interaction scheme discussed in Section~\ref{sect:binary}.

By repeating the same procedure in the case of the cost functional~\eqref{eq:J.vd} we determine the following control:
\begin{equation}
	u=\frac{\Delta{t}}{\nu+\Delta{t}^2}(v_d-v)-\frac{\Delta{t}^2}{\nu+\Delta{t}^2}I(v,\,w;\,\rho)
	\label{eq:u.vd}
\end{equation}
which finally gives rise to the feedback-constrained binary interaction scheme
\begin{align}
	\begin{aligned}[c]
		v' &= v+\dfrac{\nu\Delta{t}}{\nu+\Delta{t}^2}I(v,\,w;\,\rho)+\frac{\Delta{t}^2}{\nu+\Delta{t}^2}(v_d-v) \\
		w' &= w.
	\end{aligned}
	\label{eq:binary.u.vd}
\end{align}
Also in this case, the restriction $0<\Delta{t}\leq 1$ guarantees that $v'\in [0,\,1]$ for all $v,\,w\in [0,\,1]$ and all $\nu>0$.

\section{Boltzmann-type description}
\label{sect:boltzmann}
The constrained binary interaction rules~\eqref{eq:binary.u.variance},~\eqref{eq:binary.u.vd} can be fruitfully encoded in a Boltzmann-type statistical description of the system, which is suitable to depict the aggregate dynamics. To this end we introduce the distribution function $f=f(t,\,v):\R_+\times [0,\,1]\to\R_+$ such that $f(t,\,v)dv$ is the fraction of vehicles travelling with speed comprised between $v$ and $v+dv$ at time $t$. Under a given microscopic binary interaction scheme, the time evolution of $f$ is ruled by the following Boltzmann-type equation (cf.~\cite{pareschi2013BOOK}):
\begin{multline*}
	\frac{d}{dt}\int_0^1\varphi(v)f(t,\,v)\,dv \\
	=\frac{\rho}{2}\int_0^1\int_0^1\Bigl(\varphi(v')+\varphi(w')-\varphi(v)-\varphi(w)\Bigr) \\
		\times f(t,\,v)f(t,\,w)\,dv\,dw,
\end{multline*}
that here we have written in weak form for a test function $\varphi:[0,\,1]\to\R$. Taking $\varphi\equiv 1$ we notice that the equation implies
$$ \frac{d}{dt}\int_0^1 f(t,\,v)\,dv=0, $$
thus if $f$ is chosen to be a probability distribution in $v$ at the initial time $t=0$ it will be so at every successive time $t>0$. The physical counterpart of this fact is the conservation of the mass of vehicles.

Furthermore, with specific reference to the interaction rules~\eqref{eq:binary.u.variance},~\eqref{eq:binary.u.vd}, and in particular to the fact that $w'=w$, we observe that the Boltzmann equation specialises as
\begin{multline}
	\frac{d}{dt}\int_0^1\varphi(v)f(t,\,v)\,dv \\
	=\frac{\rho}{2}\int_0^1\int_0^1(\varphi(v')-\varphi(v))f(t,\,v)f(t,\,w)\,dv\,dw,
	\label{eq:boltzmann}
\end{multline}
which can be equivalently rewritten in strong form as
\begin{equation}\label{eq:boltzmann_strong}
 \partial_t f=Q(f,\,f), 
 \end{equation}
where
$$ Q(f,\,f)(t,\,v):=\frac{\rho}{2}\left(\int_0^1\frac{1}{'J}f(t,\,'v)f(t,\,'w)\,dw-f(t,\,v)\right) $$
is the \emph{collisional operator}. As a minor change of notation, we point out that in this formulation the symbols $'v,\,'w$ denote the pre-interaction speeds while $v,\,w$ denote the post-interaction speeds. Moreover, $'J$ is the Jacobian of the transformation from the pre- to the post-interaction speeds.

It is worth stressing that, thanks to the fact that $u$ is included in the interaction rules \eqref{eq:binary.u.variance} and \eqref{eq:binary.u.vd}, the control mechanism is naturally embedded into the kinetic equation~\eqref{eq:boltzmann}.

\subsection{Large-time trends}
\label{sect:large-time}
In order to gain some insights into the large-time trend of the solution to~\eqref{eq:boltzmann}, and particularly to ascertain the impact of the binary control strategies on the aggregate behaviour of the system, we take advantage of the \emph{quasi-invariant interaction limit} introduced by~\cite{toscani2006CMS}.

The basic idea is to investigate the asymptotic regime in which the effect of each binary interaction becomes negligible but the number of interactions per unit time is considerably high. For this we set
\begin{equation}
	\Delta{t}=\varepsilon, \qquad \nu=\nu_0\varepsilon \quad (\nu_0>0),
	\label{eq:scaling}
\end{equation}
where $\varepsilon>0$ is meant to be a small parameter, and we introduce the new time scale $\tau:=\varepsilon t$, which, owing to the scaling by $\varepsilon$, is much larger than the characteristic time scale $t$ of the binary interactions. Consequently we define the scaled distribution function $g(\tau,\,v):=f(\tau/\varepsilon,\,v)$, which from~\eqref{eq:boltzmann} is readily seen to satisfy
\begin{multline}
	\frac{d}{d\tau}\int_0^1\varphi(v)g(\tau,\,v)\,dv \\
	=\frac{\rho}{2\varepsilon}\int_0^1\int_0^1(\varphi(v')-\varphi(v))g(\tau,\,v)g(\tau,\,w)\,dv\,dw.
	\label{eq:boltzmann.g}
\end{multline}
Since for $\varepsilon$ small we have $t=\tau/\varepsilon$ large, the limit $\varepsilon\to 0^+$ describes the large-time behaviour of $f$. On the other hand, by definition of $g$, the large-time behaviour of $f$ is well approximated by that of $g$.

Let us consider, as a reference for comparison, the unconstrained interaction dynamics discussed in Section~\ref{sect:binary}. Choosing $\varphi(v)=v,\,v^2$, respectively, in~\eqref{eq:boltzmann.g} and then plugging~\eqref{eq:binary} under the scaling~\eqref{eq:scaling}$_1$ we discover, in the limit $\varepsilon\to 0^+$,
\begin{align}
	\begin{aligned}[c]
		\frac{d\cV}{d\tau} &= \frac{\rho}{2}\int_0^1\int_0^1 I(v,\,w;\,\rho)g(\tau,\,v)g(\tau,\,w)\,dv\,dw, \\
		\frac{d\cE}{d\tau} &= \rho\int_0^1\int_0^1 vI(v,\,w;\,\rho)g(\tau,\,v)g(\tau,\,w)\,dv\,dw,
	\end{aligned}
	\label{eq:dVE/dtau}
\end{align}
where
$$ \cV(\tau):=\int_0^1 vg(\tau,\,v)\,dv, \quad \cE(\tau):=\int_0^1 v^2g(\tau,\,v)\,dv $$
are the mean speed and energy of the system. Performing the same calculations with the constrained interaction rules~\eqref{eq:binary.u.variance} and the scaling~\eqref{eq:scaling}$_{1-2}$ and denoting by $V(\tau)$, $E(\tau)$ the corresponding new mean speed and energy of the system we find
\begin{equation}
	\frac{dV}{d\tau}=\frac{d\cV}{d\tau}, \qquad
	\frac{dE}{d\tau}=\frac{d\cE}{d\tau}-\frac{\rho}{\nu_0}\left(E-V^2\right)\leq\frac{d\cE}{d\tau},
	\label{eq:cal_VE}
\end{equation}
the inequality in the second equation being due to that $E-V^2\geq 0$ because this expression is the variance of the distribution $g$. If we assume that the initial speed distribution is the same in the two cases, so that $V(0)=\cV(0)$ and $E(0)=\cE(0)$, we further obtain
$$ V(\tau)=\cV(\tau), \quad E(\tau)\leq\cE(\tau), \qquad \text{for all }\tau\geq 0, $$
whence
\begin{align}
\begin{aligned}[b]
	E(\tau)-V^2(\tau) &= E(\tau)-\cV^2(\tau) \\
	& \leq\cE(\tau)-\cV^2(\tau), \qquad \text{for all }\tau\geq 0,
	\label{eq:variance_bound_w}
\end{aligned}
\end{align}
which shows that the \emph{binary} control strategy~\eqref{eq:u.variance} succeeds in reducing \emph{globally} the speed variance in the traffic flow, viz. in mitigating the component of the collective road risk linked to the differences in the speed of the vehicles. Interestingly, this happens without affecting the natural mean speed of the flow.

Similar arguments can be repeated for the binary control strategy~\eqref{eq:u.vd}, for which we obtain:
$$ \frac{dV}{d\tau}=\dfrac{d\cV}{d\tau}+\frac{\rho}{2\nu_0}\left(v_d-V\right), \quad
	\frac{dE}{d\tau}=\frac{d\cE}{d\tau}-\frac{\rho}{\nu_0}\left(E-v_dV\right). $$
For $\tau\to+\infty$, using the bounds $\abs{\frac{d\cV}{d\tau}}\leq\frac{\rho}{2}$ and $\abs{\frac{d\cE}{d\tau}}\leq\rho$ deducible from~\eqref{eq:dVE/dtau}, we can estimate $\vert V-v_d\vert\leq\nu_0$ and $\vert E-v_d^2\vert\leq\nu_0(v_d+1)$, whence $E-V^2=O(\nu_0)$. Thus also strategy~\eqref{eq:u.vd} operates so as to reduce the global speed variance of the car flow, being however more coercive than strategy~\eqref{eq:u.variance}. In fact, depending on the stregth $\nu_0$ of the control, it tends to force the mean speed towards $v_d$.

\section{Numerical examples: fundamental diagrams and speed variance}
\label{sect:numerics}

In this section we present numerical results concerning the constrained traffic model \eqref{eq:boltzmann_strong} presented in the previous sections. The results have been obtained by means of direct Monte Carlo methods for the Boltzmann equation under the scaling~\eqref{eq:scaling}, see~\cite{herty2005CMAM,pareschi2001ESAIMP,pareschi2013BOOK} for details on the numerical methods. Each test of the present section has been performed setting $\gamma=1$ in~\eqref{eq:P} and $\varepsilon=10^{-2}$ in~\eqref{eq:scaling}.

\begin{figure}
\centering
\includegraphics[scale=0.236]{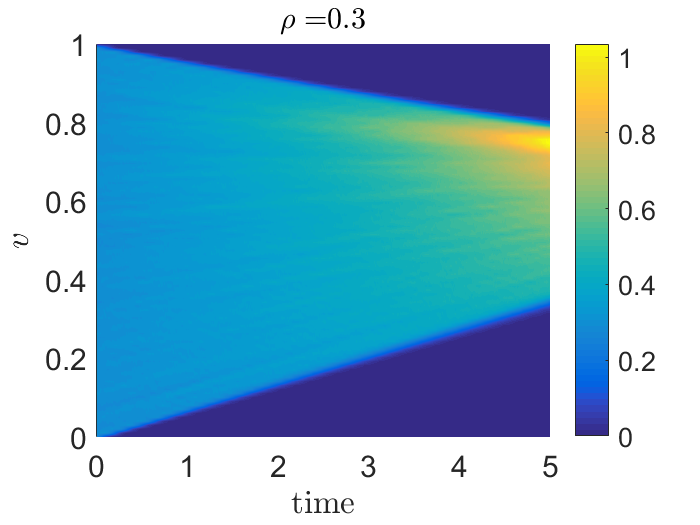}
\includegraphics[scale=0.236]{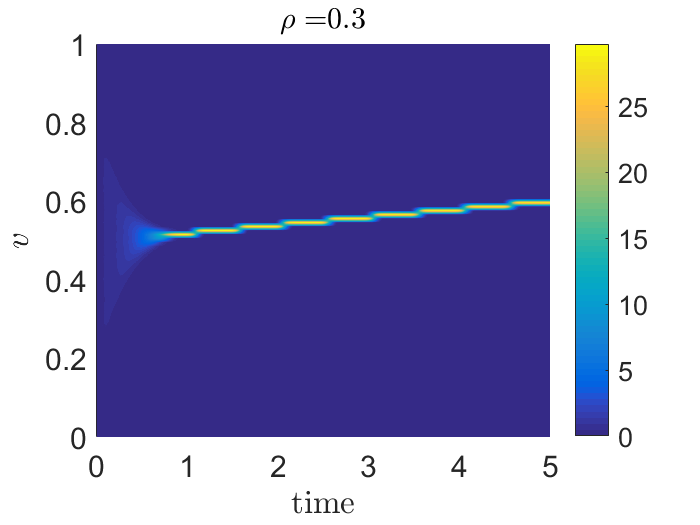} \\
\includegraphics[scale=0.236]{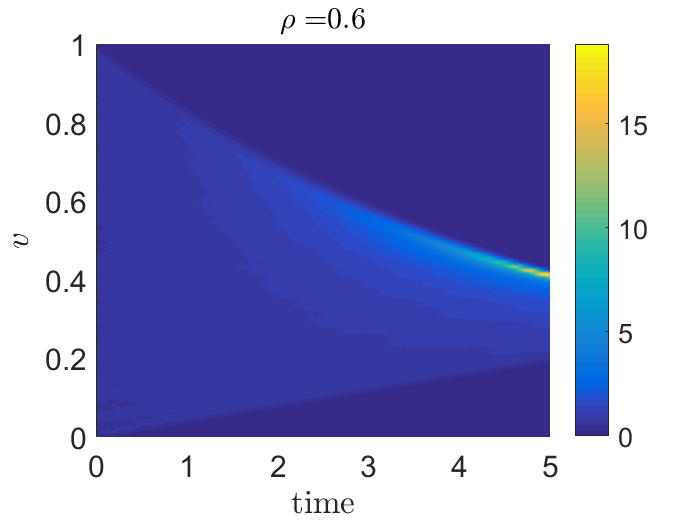}
\includegraphics[scale=0.236]{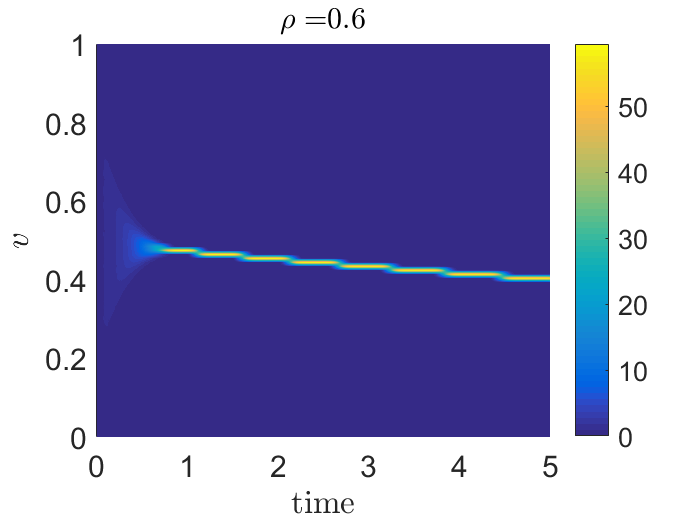}
\caption{\emph{Binary variance control}: Contours of the kinetic distribution function in the time interval $[0,\,5]$ under the unconstrained dynamics~\eqref{eq:binary} (left column) and the constrained dynamics~\eqref{eq:binary.u.variance} (right column) for the two values of the traffic density $\rho=0.3$ (first row) and $\rho = 0.6$ (second row).}
\label{fig:fevo_bin}
\end{figure}

\begin{figure}
\centering
\includegraphics[scale=0.28]{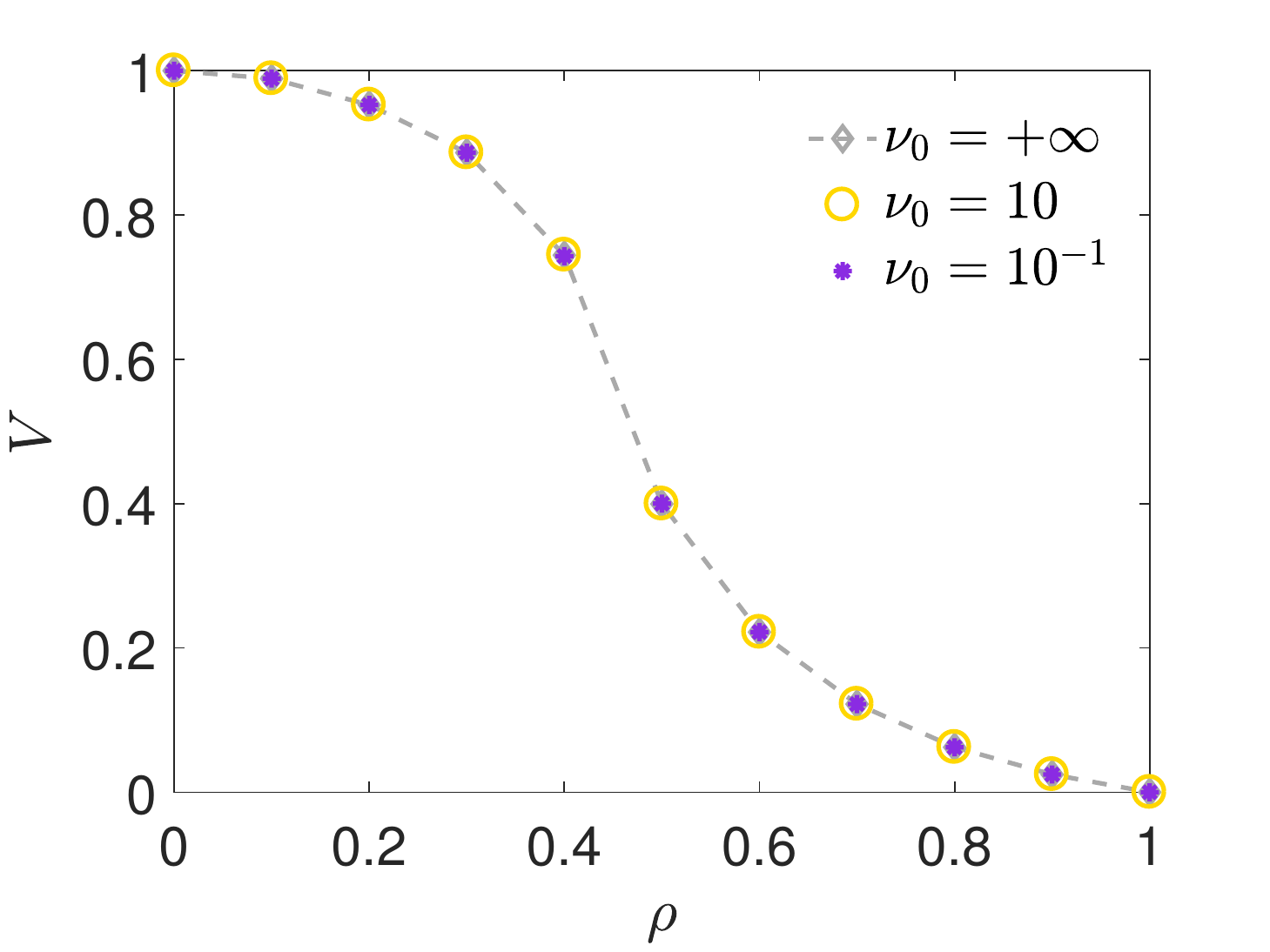}
\includegraphics[scale=0.28]{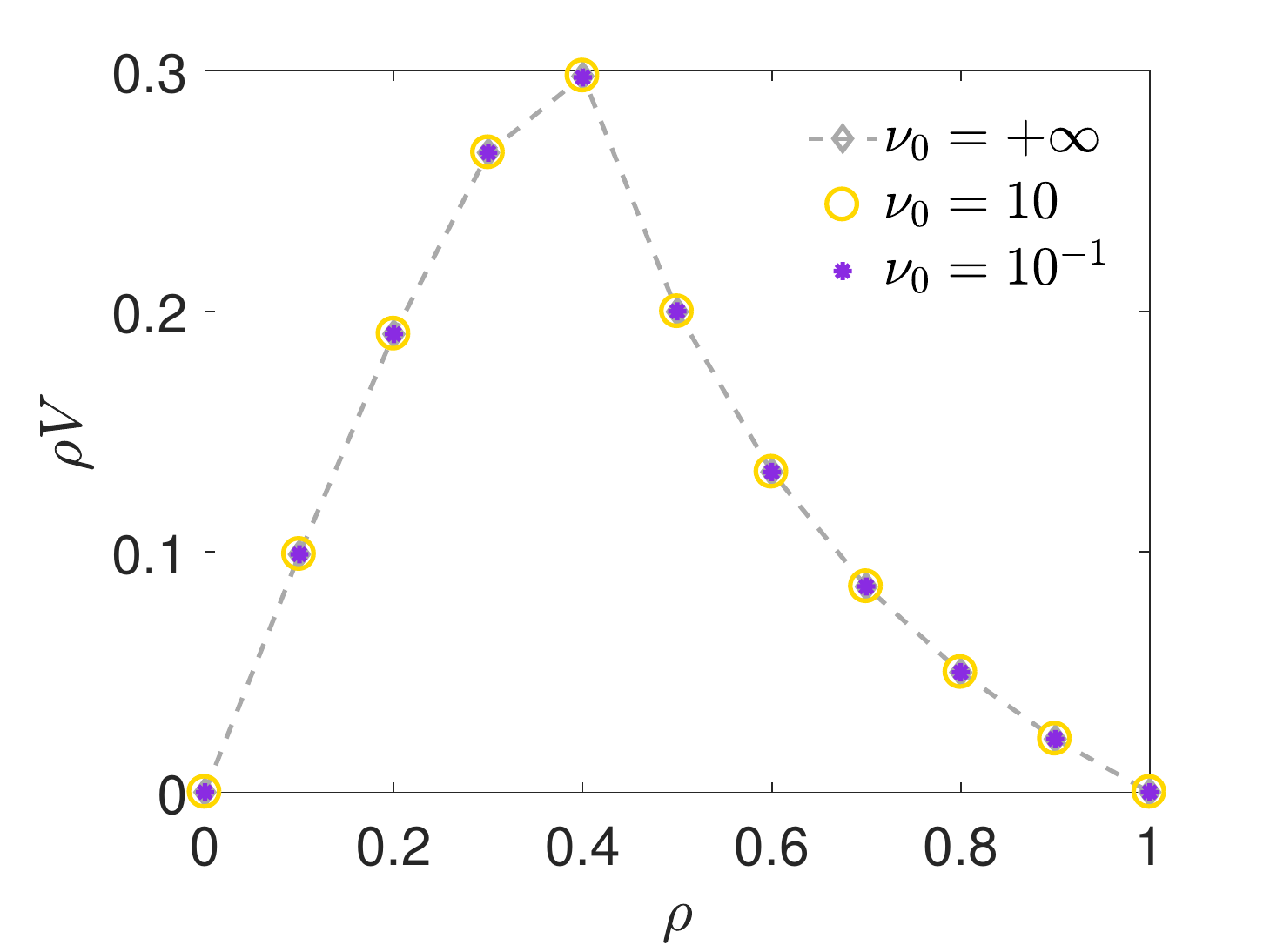}
\caption{\textit{Binary variance control}: Fundamental diagrams of traffic obtained asymptotically with the unconstrained dynamics~\eqref{eq:binary}, corresponding to $\nu_0=+\infty$, and the constrained dynamics~\eqref{eq:binary.u.variance} for two finite values of $\nu_0$. For each $\rho\in [0,\,1]$ the model has been integrated towards the steady state over the time interval $[0,\,100]$.}
\label{fig:controlw_diag}
\end{figure}

\begin{figure}
\centering
\includegraphics[scale=0.28]{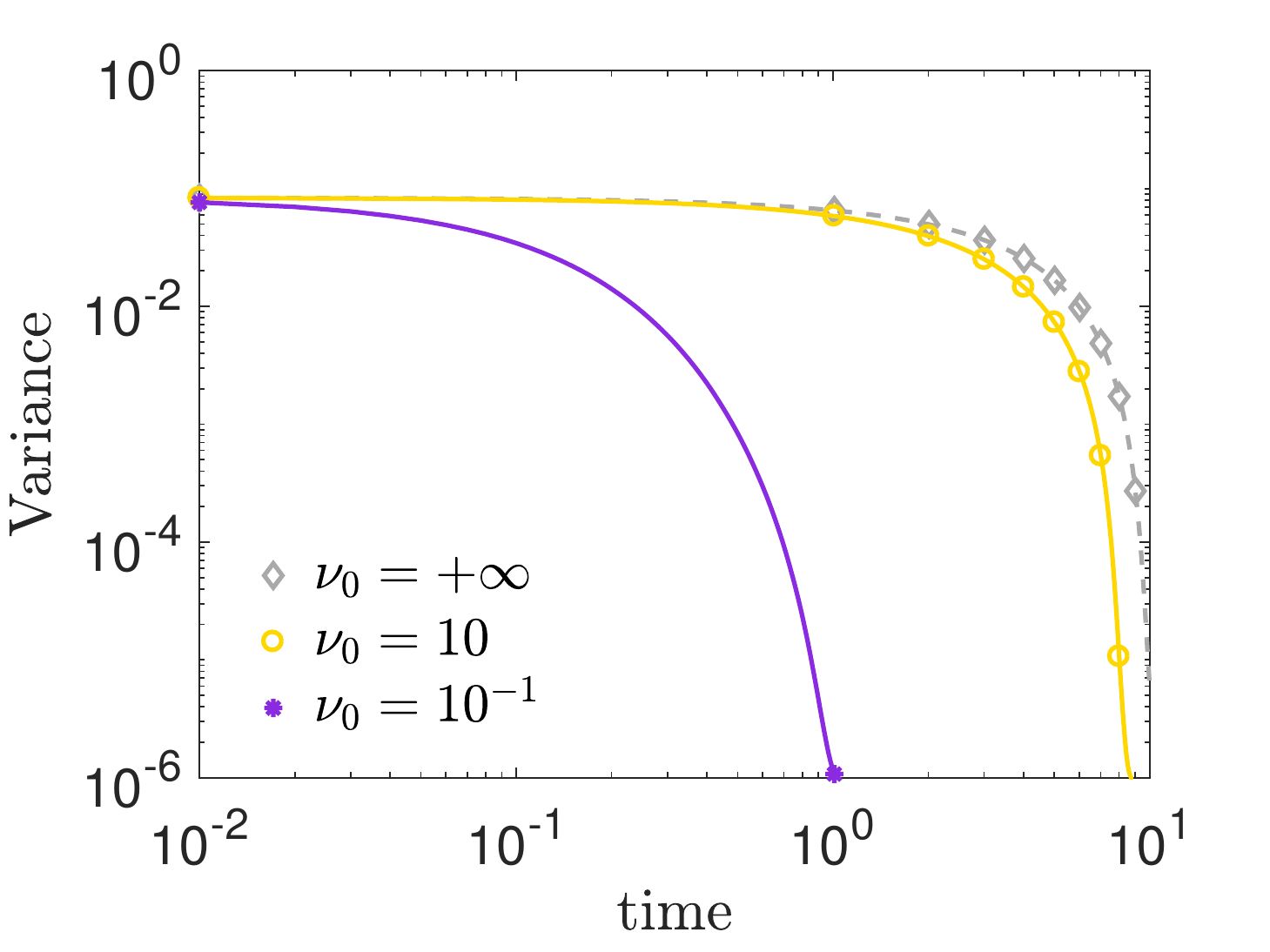}
\includegraphics[scale=0.28]{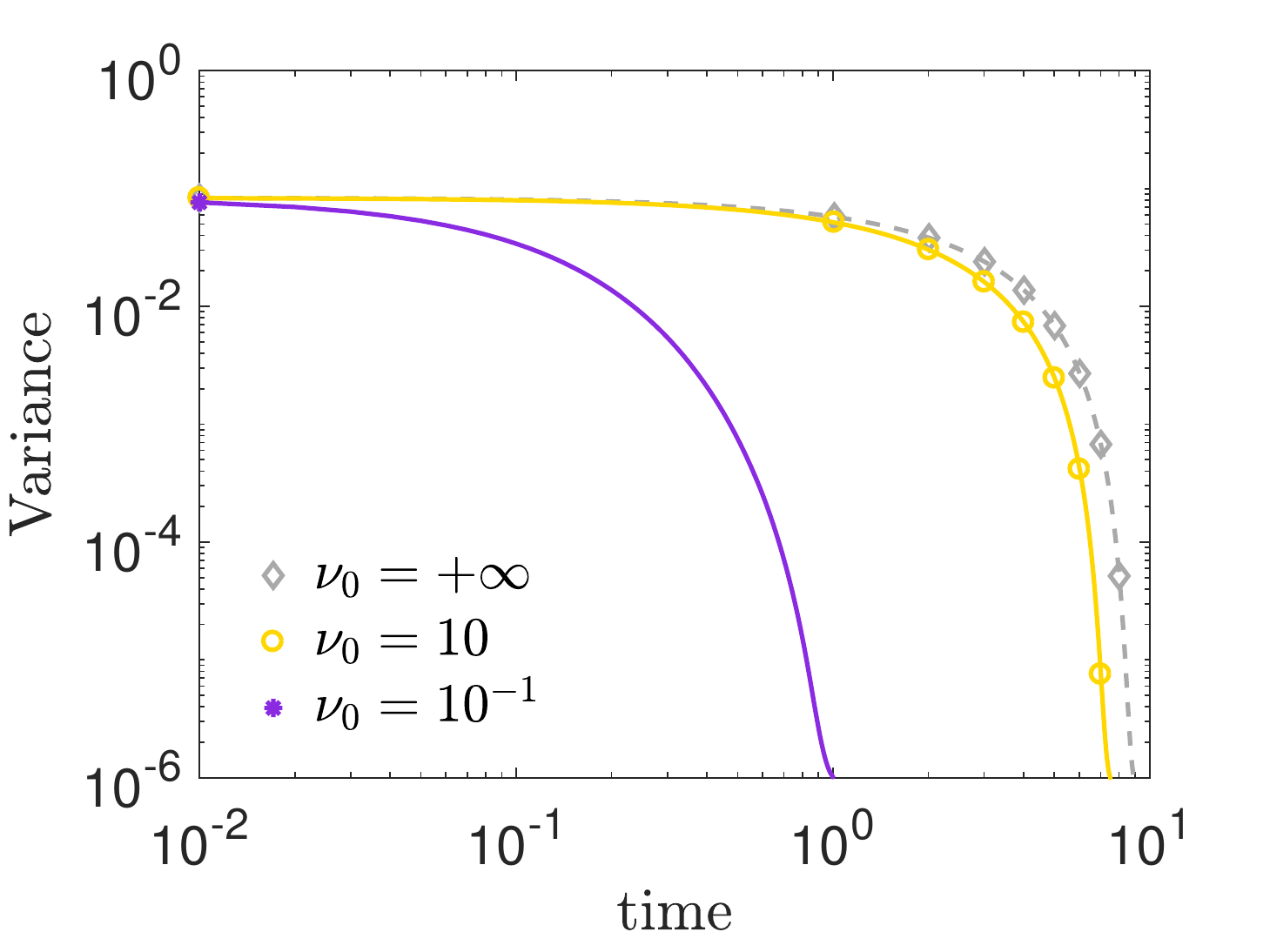}
\caption{\textit{Binary variance control}: Time evolution of the variance $E(\tau)-V^2(\tau)$, $\tau\in[0,\,10]$, for several penalisation coefficients and for $\rho=0.3$ (left) and $\rho=0.6$ (right).}
\label{fig:variance_w}
\end{figure}

\subsection{Binary variance control~\eqref{eq:J.variance},~\eqref{eq:u.variance},~\eqref{eq:binary.u.variance}}
In Figure~\ref{fig:fevo_bin} we show the contours of the kinetic distribution function in the time frame $[0,\,5]$ obtained in the unconstrained case~\eqref{eq:binary} (left column) and under the action of the binary control~\eqref{eq:u.variance} (right column) for two different values of the traffic density: $\rho=0.3$, representative of a free traffic regime, and $\rho=0.6$, representative of a congested traffic regime. It is apparent that the action of the binary control reduces immediately the variance of the speed distribution. On the other hand, from the fundamental diagrams of traffic displayed in Figure~\ref{fig:controlw_diag} we see that the control does not affect either the mean speed or the macroscopic flux of the flow of vehicles, as it has been anticipated theoretically in Section~\ref{sect:large-time}. Finally, in Figure~\ref{fig:variance_w} we show the time evolution of the variance of the speed distribution under unconstrained and constrained binary interactions and, in particular, we consider in the latter case two different values of the penalisation parameter $\nu_0$ in~\eqref{eq:scaling}: $\nu_0=10^{-1}$ (weak penalisation, strong control) and $\nu_0=10$ (strong penalisation, weak control). Consistently with the theoretical predictions, cf.~\eqref{eq:variance_bound_w}, we observe that at each time step the variance of the constrained model is bounded from above by that of the unconstrained model.

\begin{figure}
\centering
\includegraphics[scale=0.236]{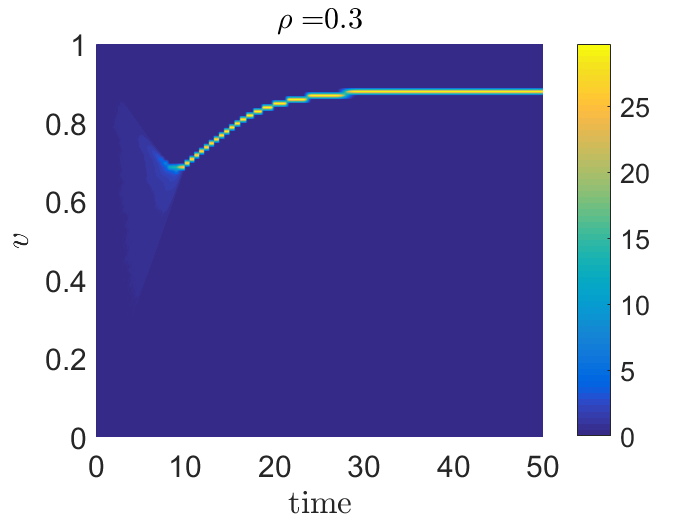}
\includegraphics[scale=0.236]{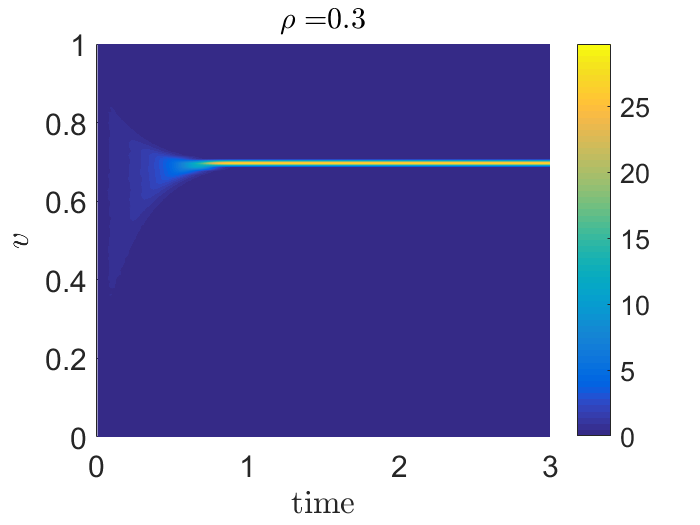} \\
\includegraphics[scale=0.236]{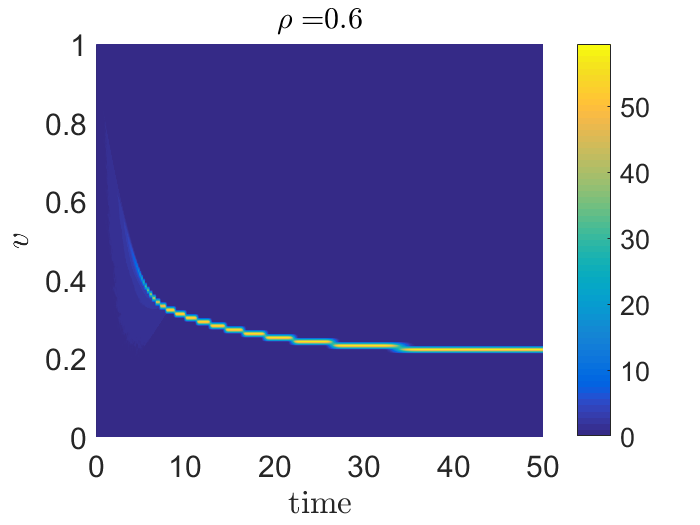}
\includegraphics[scale=0.236]{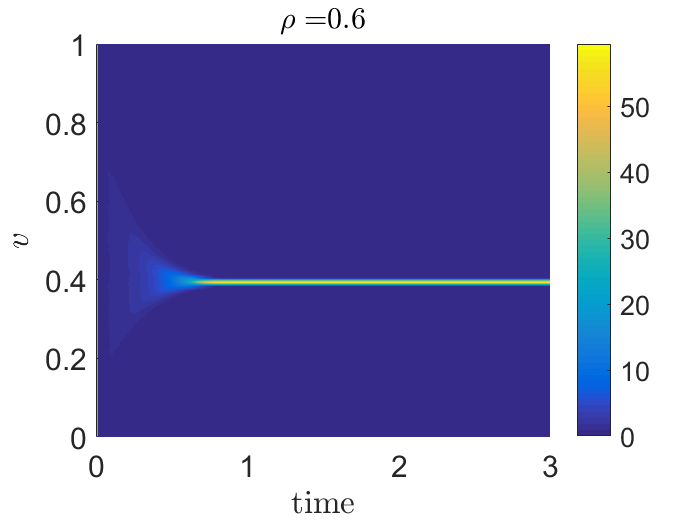}
\caption{\textit{Desired speed control}: Contours of the kinetic distribution function under the unconstrained dynamics~\eqref{eq:binary} in the time interval $[0,\,50]$ (left column) and the constrained dynamics~\eqref{eq:binary.u.vd} in the time interval $[0,\,3]$ (right column) for the two values of the traffic density $\rho=0.3$ (first row) $\rho=0.6$ (second row).}
\label{fig:fevo}
\end{figure}

\begin{figure}
\centering
\includegraphics[scale=0.28]{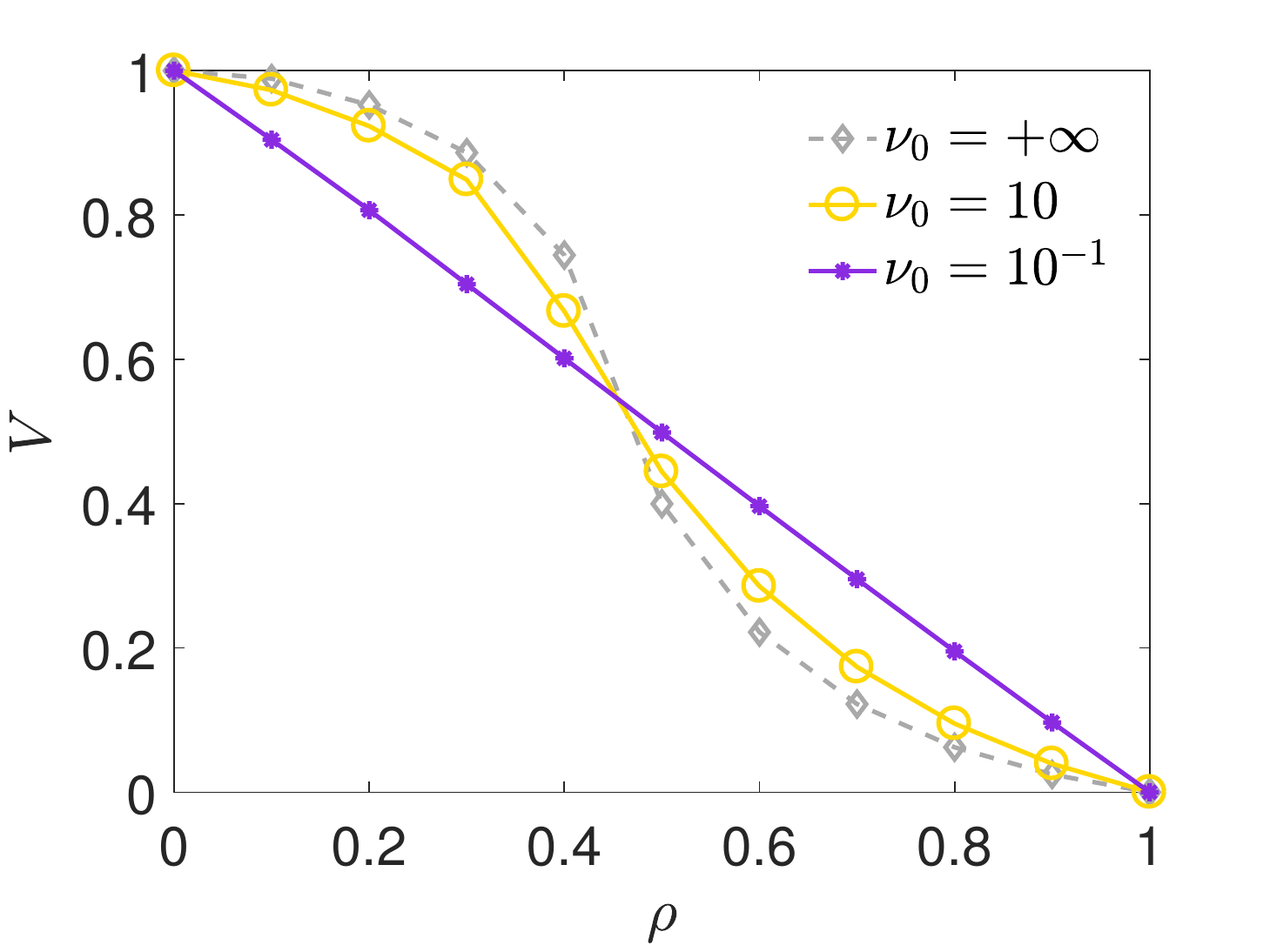}
\includegraphics[scale=0.28]{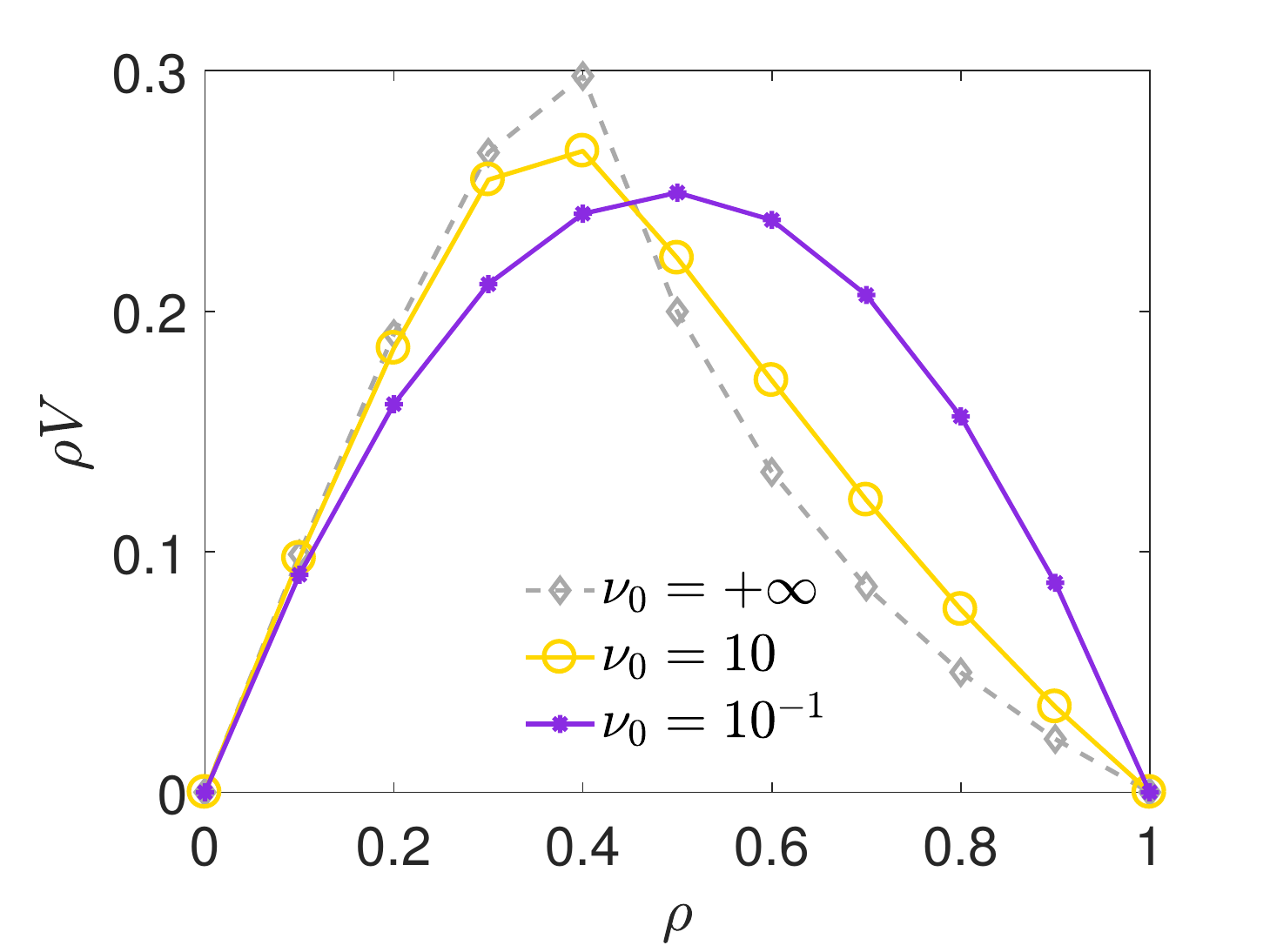}
\caption{\textit{Desired speed control}: Fundamental diagrams of traffic obtained asymptotically with the unconstrained~\eqref{eq:binary} and the constrained~\eqref{eq:binary.u.vd} dynamics for two values of the penalisation coefficient and $v_d$ like in~\eqref{eq:vd_def}. For each $\rho\in [0,\,1]$ the model has been integrated towards the steady state in the time frame $[0,\,100]$. Consistently with the theoretical findings of Section~\ref{sect:large-time}, the weaker the penalisation the closer the fundamental diagrams to that forced by the choice of $v_d(\rho)$.}
\label{fig:controlvd_diag}
\end{figure}

\begin{figure}
\centering
\includegraphics[scale=0.28]{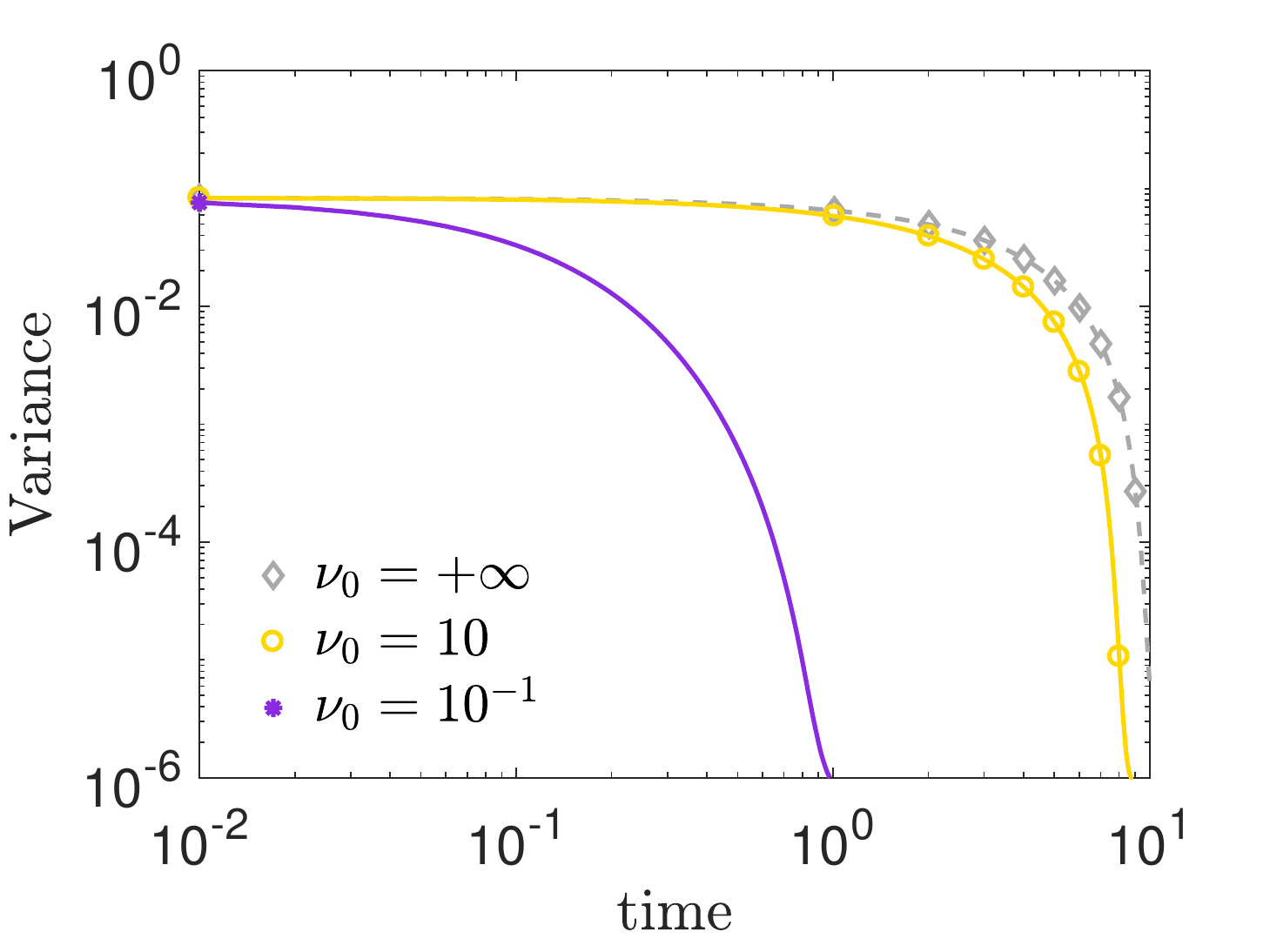}
\includegraphics[scale=0.28]{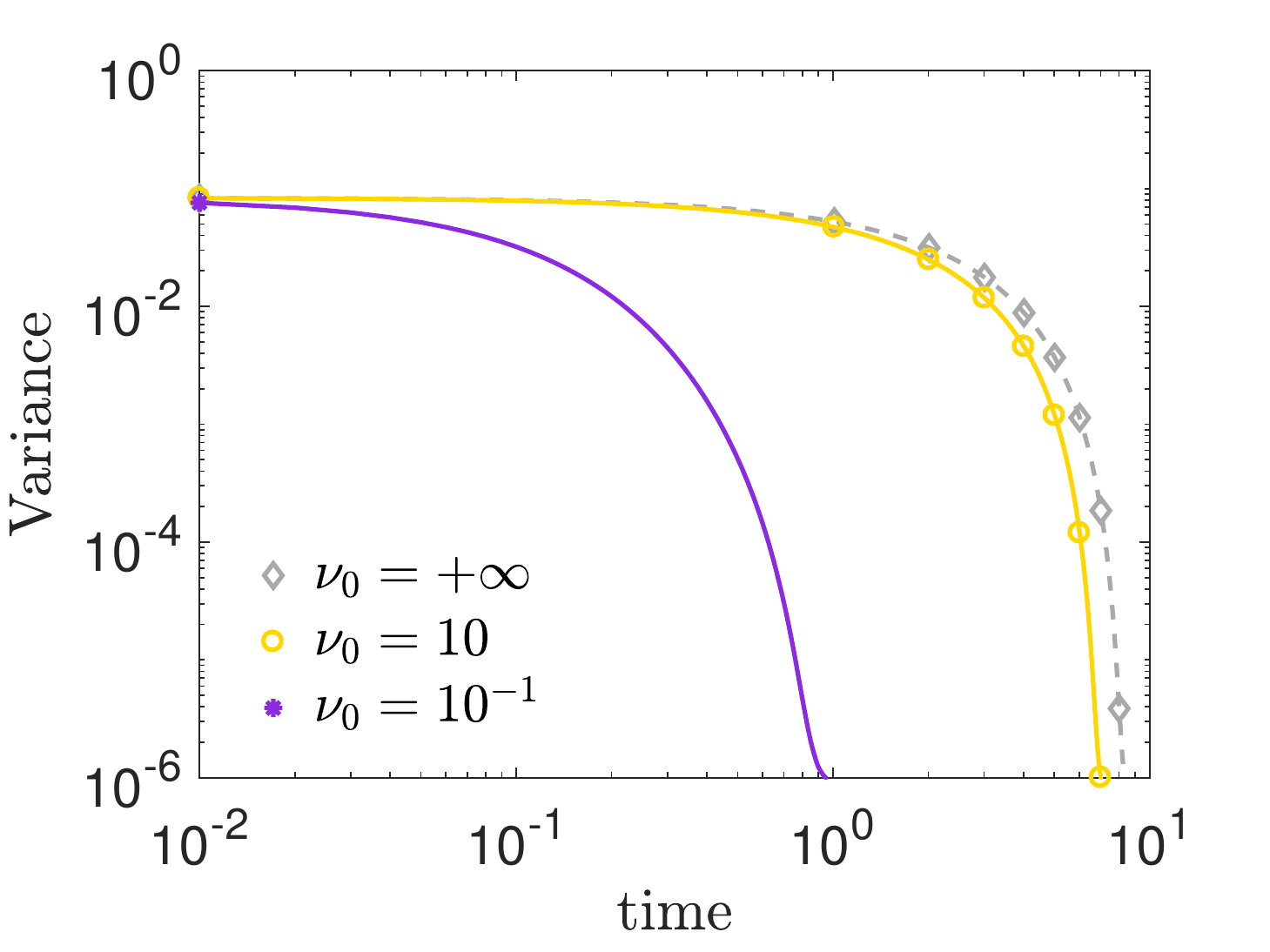}
\caption{\textit{Desired speed control}: Time evolution of the variance $E(\tau)-V^2(\tau)$, $\tau\in[0,\,10]$, for several penalisation coefficients and for $\rho=0.3$ (left) and $\rho=0.6$ (right).}
\label{fig:variance_vd}
\end{figure}

\subsection{Desired speed control~\eqref{eq:J.vd},~\eqref{eq:u.vd},~\eqref{eq:binary.u.vd}}
Concerning the control by means of the desired speed, we consider in particular a density-dependent $v_d$ of the form
\begin{equation}
	v_d=v_d(\rho)=1-\rho, \qquad \rho\in [0,\,1],
	\label{eq:vd_def}
\end{equation}
mimicking the fact that the driver-assist system may tune the target speed of the vehicle taking into account the level of congestion of the road. The relationship~\eqref{eq:vd_def} is a prototypical one implying that the desired speed is a non-increasing function of the traffic density, which vanishes in bumper-to-bumper conditions ($\rho=1$).

In Figure~\ref{fig:fevo} we show the contours of the kinetic distribution function in both the unconstrained and the constrained case, cf.~\eqref{eq:binary},~\eqref{eq:binary.u.vd}, respectively, for the same values of the traffic density $\rho=0.3,\,0.6$ as before. It is evident that the speed distribution concentrates asymptotically in two different values, the constrained one being dictated by $v_d(\rho)$ as predicted theoretically in Section~\ref{sect:large-time}. As we see from Figure~\ref{fig:controlvd_diag}, this implies that in principle such a control strategy allows one to force the fundamental diagrams of traffic to adapt to $\rho\mapsto v_d(\rho)$ (mean speed) and to $\rho\mapsto\rho v_d(\rho)$ (macroscopic flux). In particular, the choice~\eqref{eq:vd_def} of $v_d$ induces a mean speed and a macroscopic flux which are lower than the unconstrained ones in the free traffic regime (low $\rho$) but higher in the congested traffic regime (high $\rho$) while still reducing the global speed variance, hence the related road risk, at each time step, cf. Figure~\ref{fig:variance_vd}.

\section{Conclusion}
\label{sect:conclusion}
In this paper we have described a mathematical approach to control problems in kinetic traffic modelling, with particular reference to road risk mitigation issues, whose hallmarks can be summarised as follows:
\begin{inparaenum}[$(i)$]
\item the control method is based on the MPC strategy, which assumes that drivers determine their best actions by minimising a cost functional during a short and receding time horizon;
\item the time horizon is taken coincident with the duration of a single binary interaction with the leading vehicle, thereby allowing for a binary control implemented directly at the microscopic level;
\item the microscopic control problem can be solved in feedback form, i.e. the control can be expressed in terms of the microscopic states of the interacting vehicles, whereby constrained binary interaction rules can be defined explicitly;
\item the constrained binary interaction rules can be embedded in a Boltzmann-type kinetic description of the system, which allows for a statistical study of the global traffic dynamics and of the collective impact of the microscopic control strategies.
\end{inparaenum}

Starting from the consideration that the differences in the speeds of the vehicles are reported as one of the major road risk factors, we have constructed two possible control strategies for the reduction of the speed variance in the stream of vehicles. One of them does not change the fundamental diagram at the macroscopic level whereas the other drives the global mean speed towards a congestion-dependent desired speed. Both strategies have proved effective in reducing the global statistical dispersion of the speeds of the vehicles, hence potentially in mitigating the road risk component linked to the speed variance. In this preliminary approach we assumed that all vehicles are subject to the action of the control. Further developments towards more realistic scenarios may include instead sparse control strategies. 

In our view, the proposed approach can provide a sound theoretical framework to model, analyse and simulate driver-assist car technologies from a genuine multiscale perspective, with useful implications also for the traffic governance.

\begin{ack}
M.Z. acknowledges support from ``Compagnia di San Paolo'' (Torino, Italy).
\end{ack}

\bibliography{TaZm-traffic_control}
\end{document}